\documentclass[reqno]{amsart}
\usepackage[url=false,giveninits=true,maxnames=5,backend=biber]{biblatex}
\makeatletter
\bibliography{sources}
\renewcommand\part{%
	\if@noskipsec \leavevmode \fi
	\par
	\addvspace{4ex}%
	\@afterindentfalse
	\secdef\@part\@spart}

\def\@part[#1]#2{%
	\ifnum \c@secnumdepth >\m@ne
	\refstepcounter{part}%
	\addcontentsline{toc}{part}{\thepart\hspace{1em}#1}%
	\else
	\addcontentsline{toc}{part}{#1}%
	\fi
	{\parindent \z@ \raggedright
		\interlinepenalty \@M
		\normalfont
		\ifnum \c@secnumdepth >\m@ne
		\Large\bfseries \partname\nobreakspace\thepart
		\par\nobreak
		\fi
		\huge \bfseries #2%
		\par}%
	\nobreak
	\vskip 3ex
	\@afterheading}
\def\@spart#1{%
	{\parindent \z@ \raggedright
		\interlinepenalty \@M
		\normalfont
		\huge \bfseries #1\par}%
	\nobreak
	\vskip 3ex
	\@afterheading}
\makeatother
 
\usepackage{color}
\usepackage[dvipsnames]{xcolor}

\usepackage{ifpdf}
\ifpdf 
    \usepackage[pdftex]{graphicx}   % to include graphics
    \usepackage[pdftex,     % sets up hyperref to use pdftex driver
            plainpages=false,   % allows page i and 1 to exist in the same document
            breaklinks=true,    % link texts can be broken at the end of line
            colorlinks=true,
            linkcolor=red,
            citecolor=green,
            pdftitle=My Document
            pdfauthor=My Good Self
           ]{hyperref} 
\else 
    \usepackage{graphicx}       % to include graphics
    \usepackage{hyperref}       % to simplify the use of \href
\fi 

\usepackage{graphicx}
\usepackage{aurical}
\usepackage{amsfonts,amsmath}
\usepackage{amssymb}
\usepackage{verbatim}
\usepackage{amsopn}
\usepackage[english]{babel}
\usepackage{amsthm}
\usepackage{enumerate}
\usepackage{mathrsfs}
\usepackage{enumitem}
\usepackage{mathtools}
\usepackage{esint}
\usepackage{bbm}
\usepackage{marginnote}
\usepackage[marginparwidth=2cm]{geometry}
\usepackage[normalem]{ulem}
\usepackage{glossaries}
\usepackage{glossary-mcols}
\usepackage{subcaption}
\usepackage{csquotes}
\usepackage[pagewise]{lineno}%\linenumbers

\date{\today}

\theoremstyle{definition} \newtheorem{definition}{Definition}[section]
\theoremstyle{definition} \newtheorem{remark}[definition]{Remark}
\theoremstyle{plain} 
\theoremstyle{plain} 
\theoremstyle{plain} \newtheorem{theorem}[definition]{Theorem}
\theoremstyle{plain} 
\theoremstyle{definition} 
\theoremstyle{plain} 
\theoremstyle{definition} 
\theoremstyle{definition}

\DeclareMathOperator{\clos}{clos}

\newcommand{\R}{\mathbb{R}}

\newcommand{\1}{\mathbbm 1}

\newcommand{\rest}{\llcorner}

\numberwithin{equation}{section} %numera le formule in accordo con le sezioni
\theoremstyle{plain} \newtheorem*{theorem*}{Theorem}
\theoremstyle{plain} 
\theoremstyle{plain} \newtheorem*{mthm*}{Main Theorem}
\theoremstyle{plain} \newtheorem*{conjecture*}{Conjecture}
\theoremstyle{plain} 
\theoremstyle{plain} \newtheorem*{problem*}{Problem}
\title{Some counterexamples concerning a model of eradication of invasive species}

\subjclass[2020]{35K57,49K15,92-10}

\author{Stefano Bianchini}
\address{S. Bianchini: S.I.S.S.A., via Bonomea 265, 34136 Trieste, Italy}
\email{bianchin@sissa.it}

\author{Chiara Trifone}
\address{S. Trifone: S.I.S.S.A., via Bonomea 265, 34136 Trieste, Italy}
\email{ctrifone@sissa.it}

\thanks{Funding: this work has been partially supported by the GNAMPA-INDAM Grant "Problemi di controllo e regolarit\'a per evoluzione di insiemi"}

\keywords{Reaction-diffusion PDE, traveling profiles, optimal control}

\begin{document}

\maketitle

\begin{abstract}
In this note we study the regularity properties of the \emph{Effort Function $E(\beta)$}. This can be defined as the minimal $L^1$-norm of the control $\tilde \alpha(x) \geq 0$ so that the PDE
\begin{equation*}
U_t = \Delta U + f(U) - \tilde \alpha(t,x) U, \quad U \in [0,1], %\ \tilde \alpha \geq 0.
\end{equation*}
admits a traveling profile $U(x - \beta t)$ connecting $0$ to $1$. The above PDE models the evolution of an invasive species, where the control $\tilde \alpha(t,x)$ describes the percentage amount of species removal.

In the literature some conditions are assumed on the properties of $f(U)$ and (consequently) on the effort function $E(\beta)$, the main ones being that the traveling profile has a simple structure and that the effort function can be used as an anisotropic perimeter to generated a cost function over sets of finite perimeter.

This note gives some examples showing that in the interesting cases, i.e. when the PDE has a non trivial dynamics, these assumptions are in general not satisfied. We also show that $E(\beta)$ is in general only $C^1$.
\end{abstract}
%
%\keywords{}

\maketitle

%\begin{center}
%\emph{Dedicated to }
%\end{center}

\begin{center}
Preprint SISSA 16/2026/MATE
\end{center}

%\tableofcontents
%

%\newpage

%\input{senza il .tex}

\section{Introduction}
\label{S:intro}

In \cite{BressanChiri22} the following parabolic reaction diffusion equation with control is considered:
\begin{equation}
\label{Equa:parab_contr_1}
U_t = \Delta U + f(U) - \tilde \alpha(t,x) U, \quad U \in [0,1], \ \tilde \alpha \geq 0.
\end{equation}
The PDE with control \eqref{Equa:parab_contr_1} is motivated by a model of pest eradication, which fits into the large subject of optimal harvesting problems and control of parabolic equations \cite{anitacapasso2019,Ania2019GlobalEF,doi:10.1137/16M1061886,coclitegaraspino2018,coron2007,seirin2013,faucris.242204001,Lenhart2001OPTIMALCO}. The control parameter $\tilde \alpha(t,x)$ describes the percentage of invasive species removed at the point $(t,x)$.

\begin{figure}
\resizebox{\textwidth}{!}{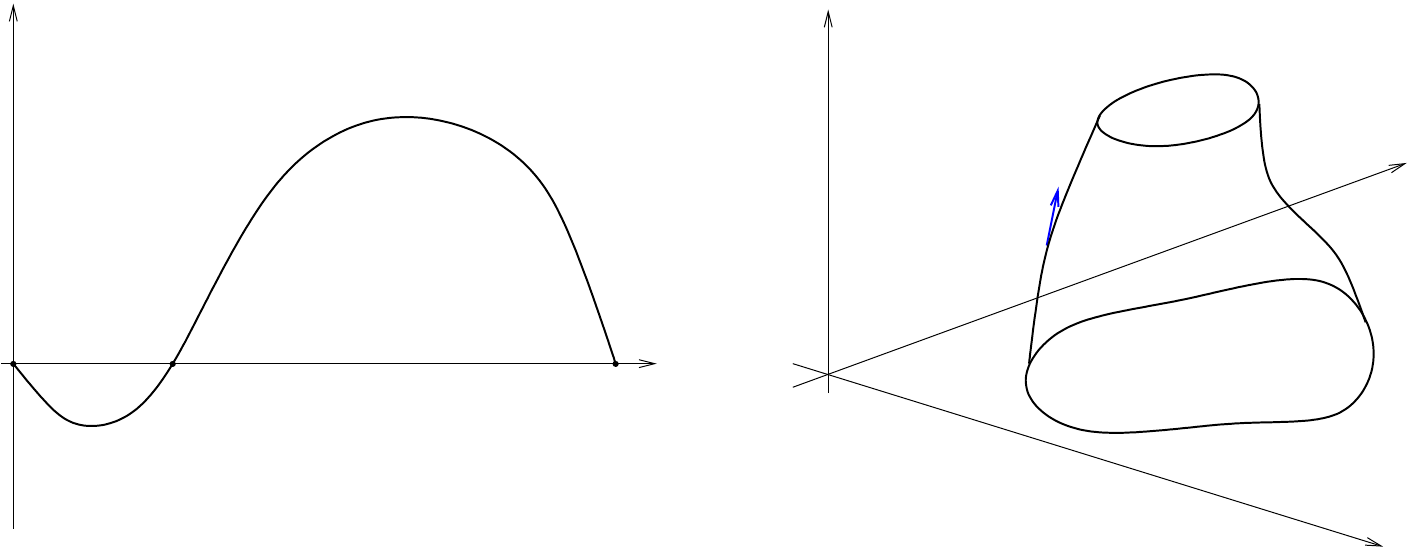}
\caption{The plot of the source $f(U)$ and the corresponding set evolution problem.}
\label{Fig:fU_Gammasu}
\end{figure}

The interesting case is when nonlinear source $f(U)$ is \emph{bistable} (see Fig. \ref{Fig:fU_Gammasu}):
\begin{enumerate}
\item $f : [0,1] \mapsto \mathbb R$ is smooth,
\item it has exactly $3$ zeros for $U = 0,U^*,1$, $0 < U^* < 1$ with $f'(0),f'(1) < 0$ and $f'(U^*) > 0$,
\item and moreover is positively predominant, i.e.
\begin{equation*}
\int_0^1 f(U) dU > 0.
\end{equation*}
\end{enumerate}
The last condition assures that the density of invasive species $U$ is growing if the control $\tilde \alpha = 0$, so that a positive control $\tilde \alpha$ is needed if the goal is to reduce or eradicate it. It can also be shown that if $\int_0^1 f(U) dU < 0$ (negatively predominant), then the corresponding set evolution problem considered in \cite{BressanChiri23} would be trivial.

One of the main results of \cite{BressanChiri22} is the analysis of the \emph{effort function $E(\beta,f)$}, which is precisely the minimal $L^1$-norm of $\tilde \alpha$ for a traveling wave solution $U(x-\beta t)$ to \eqref{Equa:parab_contr_1}. In the same paper, it is considered the rescaled problem
\begin{equation*}
\label{Equa:rescal_par}
U^\epsilon_t = \epsilon \Delta U^\epsilon + \frac{f(U^\epsilon)}{\epsilon} - \tilde \alpha^\epsilon(t,x) U^\epsilon,
\end{equation*}
and studied the possible $\Gamma$-limit of the cost functional
\begin{equation}
\label{Equa:control_intro}
\int_0^T \phi \bigg( \int_{\mathbb R^d} \tilde \alpha^\epsilon(t,x) dx \bigg) dt,
\end{equation}
with $\phi : [0,+\infty) \to [0,+\infty)$, $\phi(0) = 0$, convex and superlinear. The authors prove that if $t \mapsto A(t)$ is a smooth map into the family of smooth sets of $\mathbb R^d$, then there exists a "recovery sequence" $U^\epsilon$, solution to \eqref{Equa:rescal_par}, such that
\begin{equation*}
\lim_{\epsilon \to 0} \int_0^T \int_{\R^d} |U^\epsilon(t,x) - \1_{A(t)}| dxdt = 0
\end{equation*}
and the limit of the cost \eqref{Equa:control_intro} is
\begin{equation*}
%\label{Equa:limit_cost}
G(A) = \int_0^T \phi \bigg( \int_{\partial A(t)} E(\beta(t,x)) \mathscr H^{d-1}(dx) \bigg) dt,
\end{equation*}
where
\begin{enumerate}
\item $A$ is the smooth set
$$
A = \bigcup_{t \in [0,T]} \{t\} \times A(t) \subset \mathbb R^{d+1},
$$
\item $\beta(t,x) = - \frac{\mathbf n_t(t,x)}{|\mathbf n_x(t,x)|}$ is its normal speed at the point $(t,x) \in \partial A$,
\item the \emph{Effort Function $E(\beta)$} is computed by
\begin{equation*}
%\label{Equa:min_cost_intro}
E(\beta) = \inf \bigg\{ \int_\mathbb R \tilde \alpha(x) dx, \tilde \alpha = \frac{U_{xx} + \beta U_x + f(U)}{U}, \ U = U(x-\beta t) \ \text{traveling wave with speed $\beta$} \bigg\}.
\end{equation*}
\end{enumerate}

The fact that the expression above is the correct $\Gamma$-limit requires as necessary condition that the functional $G$ is lower semi-continuous with respect to the $L^1$-convergence of the characteristic functions of sets. Interpreting
\begin{equation*}
\int_{\partial^* A(t)} E(\beta(t,x)) \mathscr H^{d-1}(dx)
\end{equation*}
as a non standard perimeter, a sufficient condition for l.s.c. is that $\beta \mapsto E(\beta)$ is convex, positive with a linear growth at $+\infty$, see \cite{Maggi_2012}: equivalently
\begin{enumerate}
\item \label{Point1:intro_E} $E(\beta)$ is convex, l.s.c.,
\item \label{Point2:intro_E} $E(\beta) - \beta \partial_\beta E(\beta) \geq 0$.
\end{enumerate}
These assumptions are the starting point of the analysis of \cite{BressanChiri23}, and imply that the homogeneous extension
\begin{equation*}
\mathbb R \times \mathbb R^d \ni (w_t,w_x) \mapsto |w_x| E \bigg( - \frac{w_t}{|w_x|} \bigg)
\end{equation*}
is positive, convex and l.s.c., so that standard methods of calculus of variations can be applied.

In the same paper \cite{BressanChiri22}, the authors assume a monotonicity condition on the source $f(U)$. Defining the critical curve
\begin{equation*}
P^*(U) = \sqrt{Uf(U)}, \quad U \in [U^*,1],
\end{equation*}
then the monotonicity condition implies that
\begin{equation*}
\frac{d}{dU} F(U) = \frac{d}{dU} \bigg( \frac{dP^*}{dU} + f(U) \bigg) < 0.
\end{equation*}
The reason for this assumption is that the structure of traveling waves for \eqref{Equa:parab_contr_1} is particularly simple, since the control is applied only when $U$ belongs to an interval of values $U \in (U^u(\beta),U^s(\beta))$ (see Section \ref{Ss:travelin_alpha} and Theorem \ref{Theo:regul_optP_intro} for the precise definitions). For the simplest nonlinear source
\begin{equation*}
f(U) = U (U - U^*) (1 - U).
\end{equation*}
this condition is verified for $U^*$ close to $1$. However in this case $f(U)$ is not positively predominant (which requires $U^* < 0.5$), so that it is not of practical interest.

In this note we show that in general the above assumptions on the effort function $E(\beta)$ are not satisfied. We provide 4 examples addressing the monotonicity of $F(U)$, the convexity of $E(\beta)$ and its regularity: more precisely, we show that
\begin{itemize}
\item $F(U)$ is in general non monotone: as a consequence the structure of optimal profiles is more complicated that the one in Fig. \ref{Fig:optproof_1};

\item $E(\beta)$ is not convex: we give two different examples, one based on a direct analysis of the derivative of the effort function $E(\beta)$ and the other based on the existence of an asymptote as $\beta \to +\infty$;

\item $E(\beta)$ is not smooth: we show that in general the regularity to expect is $C^{1,1}$. A more refined analysis on the behavior of the optimal profile near $U = U^*$ proves that $E(\beta)$ is only $C^1$, but this requires additional techniques.
\end{itemize}

In order to understand the examples we provide, in Section \ref{S:struct_opt} we collect the key results on traveling profiles for \eqref{Equa:parab_contr_1}. \\
In Section \ref{Ss:travelin_no_alpha} we recall some known properties of traveling profiles for the PDE \eqref{Equa:parab_contr_1}. In particular, putting the ODE for traveling profiles
\begin{equation*}
- \beta U_x = U_{xx} + f(U) - \tilde \alpha(x) U
\end{equation*}
in the phase space $(U,P=U_x)$, we introduce the unstable manifold $\Gamma_u$ of $(0,0)$ and the stable manifold $\Gamma_s$ of $(1,0)$ and show that the positively predominant condition implies that speed of traveling profiles when $\tilde \alpha = 0$ is a constant $\beta^* < 0$. \\
In Section \ref{Ss:travelin_alpha} we introduce the \emph{Effort Function $E(\beta)$} on curves of the phase space, and in Section \ref{Ss:stru_opt} we show that under the monotonicity condition the optimal profile is easily computed: it is made by an arc of the unstable manifold $\Gamma_u(U)$, and an arc of the critical curve $P^*(U) = \sqrt{U f(U)}$, and an arc of the stable manifold $\Gamma_s(U)$. \\
Finally, we recall the main results in the general case, i.e. when the monotonicity assumption fails: the properties of a minimizer and its uniqueness (Theorem \ref{Theo:regul_optP_intro}), and the regularity of the effort function $E(\beta)$ (Theorem \ref{Theo:Ebeta_intro}).

After these preliminaries, we are ready to present the 4 examples.

In Section \ref{S:monoton_cond} we show that the monotonicity assumption on $f(U)$ is not satisfied even for the cubic source when $U^* < .4$. Hence, in general the structure of optimal traveling profiles is more complicated than the one considered in Section \ref{Ss:stru_opt}. A complete analysis of their structure is done in \cite{bianchini2026travelingprofilescontrolcost}, where also Theorem \ref{Theo:regul_optP_intro} is proved.

The example of Section \ref{S:E_not_convex} considers a perturbation of the Hamiltonian system when $\beta=0$, and in this case its associated manifolds can be explicitly computed, see Equations \eqref{unstable0} and \eqref{stable0}. Similarly, for $\beta=0$, we are able to provide formulas for $\frac{\partial \Gamma_u}{\partial \beta}$, $\frac{\partial \Gamma_s}{\partial \beta}$, $\frac{\partial ^2 \Gamma_u }{\partial \beta^2} $ and $\frac{\partial ^2 \Gamma_s }{\partial \beta^2}$. In Equation \eqref{Eq:nonconvsource}, we introduce a sequence $\{f_K\}_{K \in \mathbb{N}}$ of continuous sources for which we are able to show 
\begin{itemize}
    \item $\frac{\partial ^2 \Gamma_u}{\partial \beta^2}(0,U^*)>0$ and $ \frac{\partial ^2 \Gamma_u}{\partial \beta^2}(0,U_u(0))=\frac{\partial ^2 \Gamma_u}{\partial \beta^2}(0,U^*) + \mathcal O ( \frac{1}{\sqrt{K}} )$,
    \item $\frac{\partial^2\Gamma_s}{\partial \beta^2}(0, U_s(0))=\frac{1}{8\sqrt{K(1-2U^*)}} \to 0 \text{ as } K\to +\infty$,
\end{itemize}
so that, in particular, for $K \gg 1$,
\begin{equation*}
\frac{1}{U_s(0)} \frac{\partial^2 \Gamma_s}{\partial \beta^2}(0,U_s(0)) - \frac{1}{U_u(0)}\frac{\partial^2 \Gamma_u}{\partial \beta^2}(0,U_u(0))<0.
\end{equation*}
We remark that 
\begin{equation*}
    \int_0^1 f_K(U)\, dU \to +\infty \quad \text{for} \ K\to+\infty, 
\end{equation*}
meaning that requiring that $K$ is large enough corresponds  to asking that $f$ is highly positively predominant.

In Section \ref{S:anisotr_nonconvex}, we use that \eqref{Equa:asympt} gives the asymptotic linear behavior
\begin{equation*}
     E(\beta)\sim \ln \bigg( \frac{1}{U^*} \bigg) \beta +\int_{U^*}^1 \frac{2 \sqrt{f(U)}}{U^{3/2}} dU \quad \text{as}\ \beta\to+\infty.
\end{equation*}
Such property obviously has consequences on the convexity of $E$, in particular on the assumptions \eqref{Point1:intro_E}, \eqref{Point2:intro_E} of page \pageref{Point1:intro_E}. We show first that if one assumes a priori that $E$ is convex, its asymptotic behavior automatically implies the second condition. Conversely, exploiting similar arguments and the Lipschitz continuity of the effort stated in Theorem \ref{Theo:Ebeta_intro}, we are able to exhibit a cubic source whose associated effort is not convex.

The example of Section \ref{S:exampl_C11} shows that in general the effort function $E(\beta)$ is not $C^2$. Indeed, the regularity of the optimal profile $P_\beta$, and thus of $E(\beta)$, is strictly linked to the topological features of the free set $\mathcal I(P_\beta)$: such observation is crucial to construct a source whose associated effort is not of class $C^2$. We first recall the explicit expressions of the first and second derivative of $E(\beta)$ when the dependence of the control set $\mathcal O(P_\beta)$ w.r.t. $\beta$ is sufficiently regular, in particular when the number of connected components of the free set $\mathcal I(P_\beta)$ is finite and constant. We do not proved an explicit source $f(U)$ for which the number of connected components varies: such an example can be explicitly constructed and numerically studied. \\
We consider in particular the situation where a splitting occurs in $\mathcal{I}(\overline{\beta})$, more precisely
\begin{equation*}
     \mathcal{I}(\beta) = \begin{cases}
         (U_1^-(\beta),U_1^+(\beta)) & \text{for }   \beta<\overline{\beta}, \\
         (U_2^-(\beta), U_2^+(\beta)) \cup (U_3^-(\beta), U_3^+(\beta)) & \text{for }   \beta>\overline{\beta}.
     \end{cases}
\end{equation*}
Knowing that the optimal profile solves \eqref{functeq} in each connected component of $\mathcal{I}$, equals $P^*$ at the endpoints and must satisfy the necessary condition of Theorem \ref{Theo:regul_optP_intro} below, we are able to compute the right and left limits of $\partial^2_\beta E$ at $\bar \beta$ and prove that they do not coincide. This gives that at most $E(\beta)$ is $C^{1,1}$.

\section{Existence and structure of optimal profiles}
\label{S:struct_opt}

Aim of this section is to recall some basic facts about traveling profiles for parabolic reaction diffusion equations: we refer to \cite[Chapter 4]{Fife}.

\subsection{Traveling waves with no control}
\label{Ss:travelin_no_alpha}

Consider the PDE \eqref{Equa:parab_contr_1} without control
\begin{equation}
\label{eq}
U_t = U_{xx} + f(U).
\end{equation}
As in \cite{BressanChiri22} and \cite{BressanChiri23} we make the following assumptions on the \emph{source function} $f(U)$:
\begin{enumerate}[label=\textbf{(H\arabic*)}, ref=H\arabic*]
\item\label{H1} The function $f \in C^2([0,1])$ satisfies (see Fig. \ref{Fig:fU_Gammasu})
\begin{equation*}
f(0)=f(1)=0, \ f'(0)<0, \ f'(1)<0.
\end{equation*}
and there exists a unique point
$U^* \in (0,1)$ such that $f(U^*)=0$ and $f'(U^*)>0$.
\end{enumerate}
\noindent Moreover we require that $f$ is \emph{positively predominant}, that is
\begin{equation*}
%\label{Equa:positive_prev}
\int_0^1 f(U) \, dU > 0.
\end{equation*}

A \emph{traveling wave solution with velocity
$\beta$} (or \emph{traveling profile}) for \eqref{eq} is a solution of the form
\begin{equation*}
%\label{Equa:trav_sol}
U(t,x) = U(x-\beta t)
\end{equation*}
connecting $U(-\infty) = 0$ to $U(+\infty) = 1$. If we consider the  system in the phase plane $(U,P)$
\begin{align}
\label{sys}
\begin{cases}
U' = P, \\
P' = - f(U) -\beta P,
\end{cases}
\end{align}
a traveling profile corresponds to a trajectory in the $(U,P)$-plane connecting the equilibria $(0,0)$ and $(1,0)$ with $P > 0$. %\newline

Linearizing the system at a point $(U,P)$ we obtain the Jacobian matrix
\begin{equation*}
D(U,P)=\begin{pmatrix}
0 & 1 \\
- f'(U) & - \beta
\end{pmatrix},
\end{equation*}
whose eigenvalues are
\begin{align*}
%\label{Equa:eigenv_equi}
\lambda_\pm(\beta,U) = - \frac{\beta}{2} \pm \sqrt{\frac{\beta^2}{4}-f'(U)}.
\end{align*}
%{\color{blue} MA SERVE? Note that
%\begin{align*}
%\det(DF(U,P))=f'(U), \quad \Tr(DF(U<P)) = - \beta \quad \text{and} \quad \Delta DF(U,P)=\sqrt{\beta^2-4f'(U)}.
%\end{align*}}
We can therefore distinguish the following cases with respect to the stability of the fixed points $(0,0)$, $(U^*,0)$, and $(1,0)$.

\begin{enumerate}
\item When $f'(U)<0$, the eigenvalues are real and have opposite signs for all $\beta \in \R$: in particular the equilibria $(0,0)$ and $(1,0)$ are saddles.

\item If $f'(U)>0$ and $\beta <- 2 \sqrt{f'(U)}<0 $, the eigenvalues are both real and positive: the equilibrium $(U^*,0)$ is an unstable node.

\item If $f'(U)>0$ and $- 2\sqrt{f'(U)}<\beta < 0$, $\lambda_{\pm}(\beta, U)$ are complex conjugate, both with positive real part: the equilibrium $(U^*,0)$ is an unstable focus.

\item For $\beta = 0$ the system is Hamiltonian, so that $(U^*,0)$ is a center.

\item \label{Point_5:eigen_equi} If $f'(U)>0$ and $0 < \beta < 2\sqrt{f'(U)}$, $\lambda_{\pm}(\beta, U)$ are complex conjugate, both with negative real part: the equilibrium $(U^*,0)$ is a stable focus.

\item If $f'(U)>0$ and $\beta > 2\sqrt{f'(U)}>0 $, the eigenvalues are both real and negative: the equilibrium $(U^*,0)$ is a stable node.
\end{enumerate}

Therefore, for $P>0$ the phase portrait for the system \eqref{sys} contains an unstable manifold $\Gamma_u = \Gamma_u(\beta,f,U)$ of the point $(0,0)$ and a stable manifold $\Gamma_s = \Gamma_s(\beta,f,U)$ of the point $(1,0)$. Moreover (\cite[Chapter 4.4.a]{Fife}) there exists a unique value $\beta^* = \beta^*(f)$ such that $\Gamma_u(\beta^*)=\Gamma_s(\beta^*)$: thus, we can find a unique (up to translation) traveling profile $U$ with speed $\beta^*$. %This means that, in absence of control, the solution $U(x - \beta^* t)$ keeps expanding.

We write $\Gamma_u(\beta)$ as a function $\Gamma_u(\beta,U)$ in the region $[0,U_{uf}(\beta) = U_{uf}(\beta,f)]$ where $\Gamma_u(\beta) \geq 0$, with $\Gamma_u(\beta,0) = \Gamma_u(\beta,U_{uf}(\beta)) = 0$ and strictly positive inside the interval; similarly we write $\Gamma_s(\beta) = \Gamma_s(\beta,U)$, which it is easy to see that it is a function for $\beta \geq \beta^*$.

%From the monotonicity of $\beta \mapsto \Gamma_u(\beta)$, it follows that the $U_{uf}(\beta)$ is decreasing, until it eventually becomes $U^*$: 
Point \eqref{Point_5:eigen_equi} above and the fact that for $\beta \to \infty$ it holds $\Gamma_u(\beta,U) \to 0$ yields the existence of a value $\beta^{**} = \beta^{**}(f) \geq 2\sqrt{f'(U^*)}$ such that for $U_{uf}(\beta^{**}) = U^*$. From this value on $U_{uf}(\beta) = U^*$.

\begin{remark}
\label{Rem:ether_1}
The value of $\beta^{**}$ corresponds to the "heteroclinic connection" where $\Gamma_u$ coincides with the "strongly" stable manifold defined by the eigenvalue $\lambda_-(\beta,U^*)$.
\end{remark}

\begin{figure}
\resizebox{\textwidth}{!}{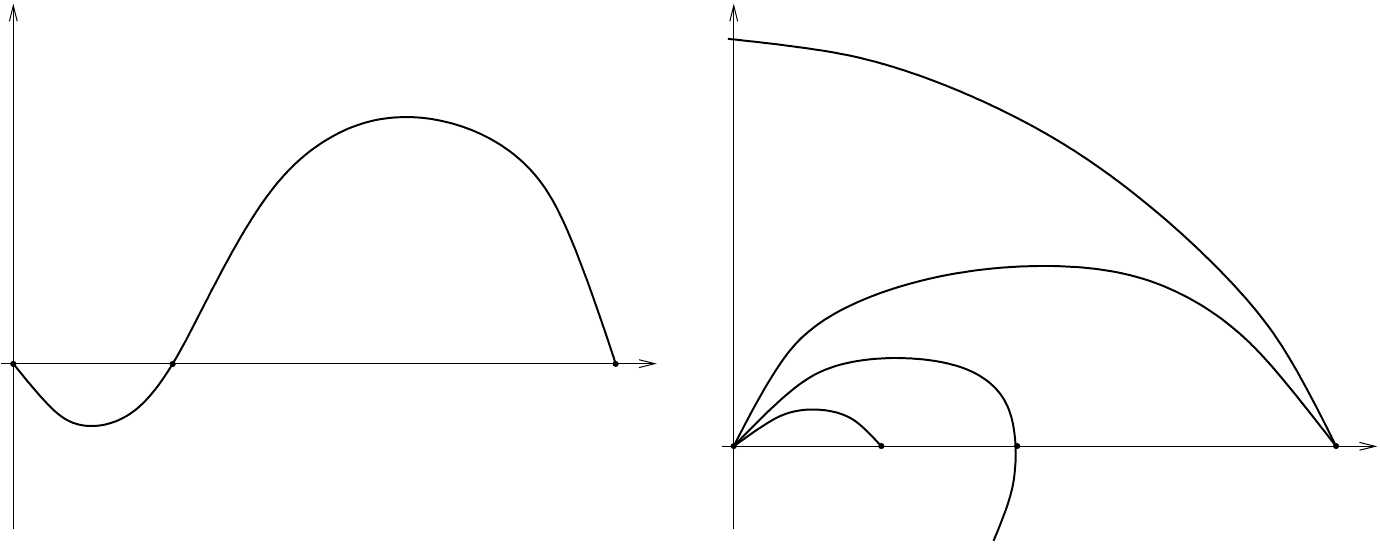}
\caption{The plot of the source $f(U)$ and of the stable and unstable manifolds for different speeds $\beta$.}
\label{Fig:fU_Gammasu_1}
\end{figure}

See Fig. \ref{Fig:fU_Gammasu_1} for a graphical representation of the manifolds $\Gamma_u(\beta),\Gamma_s(\beta)$ for different values of $\beta$. %the various cases. 

\subsection{Traveling waves with control}
\label{Ss:travelin_alpha}

Next, consider the model for PDE with control
\begin{equation}
\label{eqwithcontrol}
U_t = U_{xx} + f(U) - \tilde \alpha(t,x) U.
\end{equation}
The control $\tilde \alpha(t,x)$ is assumed to be positive, which means that our goal is to diminish $U$, and it allows to construct traveling profiles with speed $\beta > \beta^*$: for a given $\beta > \beta^*$, we look for a control $\tilde \alpha(s)$ and a function $U = U(s)$ such that $U(t,x) = U(x-\beta t), \ \tilde \alpha(t,x) = \tilde \alpha(x-\beta t)$ solve \eqref{eqwithcontrol}.

\begin{definition}[Admissible profiles]
\label{admcurves}
Given $\beta > \beta^*$, we denote by  $\mathcal A_\beta = \mathcal A_{\beta,f}$ the set of all Lipschitz curves
$$
s \mapsto \gamma(s)=(U(s), P(s)) \in (0,1) \times (0,+\infty)
$$
%for some interval $[0, \overline{s}]$ which %such that
%\begin{equation*}
%\gamma(0)=(0,0), \quad \gamma(\overline{s})=(1,0),
%\end{equation*}
%\begin{equation*}
%P(s) = U'(s) \ge 0 \quad \text{and} \quad P'(s) \ge -f(U(s))-\beta P(s) \ \forall s \in [0,\overline{s}], 
%\end{equation*}
%i.e. $\gamma$ is 
which are graphs of solutions $(U(s),P(s))$ for the planar system
\begin{align*}
%\label{syswithcontrol}
\begin{cases}
U'= P, \\
P'= -\beta P - f(U) + \tilde \alpha U,
\end{cases} \quad (U,P)(-\infty) = (0,0), \ (U,P)(+\infty) = (1,0),
\end{align*}
for some admissible control $\tilde \alpha \ge 0$.
\end{definition}

Given an admissible curve $\gamma \in \mathcal{A(\beta)}$, its associated cost
\begin{equation*}
\tilde E(\gamma) = \int_\R \tilde \alpha(x) dx = \int_{\gamma} \mathbf{v},
\end{equation*}
with
\begin{equation*}
\mathbf{v} = \left( \frac{f(U)}{UP} + \frac{\beta}{U}\right)dU + \left(\frac{1}{U}\right) dP,
\end{equation*}
i.e. a path integral along $\gamma$.

\begin{definition}[Effort function]
\label{Def:effort_funvc}
We define the \emph{effort function} $E(\beta) = E(\beta,f)$ as
\begin{equation}
\label{effort}
\begin{array}{ccccc}
E &\colon& [\beta^*, +\infty) &\to& \R^+ \\
&& \beta &\mapsto& E(\beta) = \inf \{ \tilde E(\gamma), \gamma \in \mathcal{A}_\beta\}
\end{array}
\end{equation}
\end{definition}

We can assume that $E(\beta) = +\infty$ for $\beta < \beta^*$, in the sense that there are no admissible solutions.

%\begin{remark}
%\label{Rem:control_on_UP}
Being $P > 0$ if $\gamma \in \mathcal{A}_\beta$, it follows that $s \mapsto U(s)$ is invertible in $I$ and we can express $P$ as a function of $U$, solving
\begin{equation}
\label{Equa:funct_eq_UP}
\frac{dP}{dU} = - \beta - \frac{f(U)}{P} + \alpha(U) U, \quad \alpha(U) = \frac{\tilde \alpha(s(U))}{P(U)},
\end{equation}
with
\begin{equation*}
\int_I \tilde \alpha(s) \, ds = \int_{U(I)} \frac{\tilde \alpha(s(U))}{P(U)} \, dU=\int_{U(I)} \alpha(U) \, dU.
\end{equation*}
We therefore can use the cost $\alpha(U)$ for traveling profiles when computing the cost: this will remove the degeneracy of the invariance of the problem w.r.t. translations. If $\alpha = 0$ we obtain the ODE
\begin{equation}
\label{functeq}
\frac{dP}{dU} = - \frac{f(U)}{P(U)} - \beta.
\end{equation}

\subsection{Structure of optimal profiles under the monotonicity assumption}
\label{Ss:stru_opt}

In this section we recall the results developed in \cite{BressanChiri22}, where the authors make the \emph{monotonicity assumption}
\begin{enumerate}[label=\textbf{(H\arabic*)}, ref=H\arabic*]
\setcounter{enumi}{1}
\item\label{mon}
for all $U \in [U^*,1]$, one has $(4-2\sqrt{3})f'(U) + 2U f''(U) \le 0$.
\end{enumerate}
Under this condition, it is proved that the infimum in \eqref{effort} is achieved by the smooth curve $\gamma_\beta$ obtained via the concatenation of three parts:
\begin{itemize}
\item the unstable manifold $\Gamma_u=\Gamma_u(\beta)$ until it reaches a point $(U_u(\beta), P^*(U_u(\beta)))$ on the critical curve $P^* = \sqrt{Uf(U)}$;
\item the stable manifold $\Gamma_s=\Gamma_s(\beta)$ until the point $(U_s(\beta), P^*(U_s(\beta)))$ of intersection with the graph of $P^*$;
\item the segment of the graph of $P^*$ connecting the previous two curves.
\end{itemize}
Recall that $\Gamma_u(\beta)$ and $\Gamma_s(\beta)$ are the unstable manifold of $(0,0)$ and the stable manifold of $(1,0)$. Due to the monotonicity assumption, there are no other intersections of $\Gamma_u$ or $\Gamma_s$ with $P^*$ beyond $U_u,U_s$.

\begin{figure}
\centering
\resizebox{.75\textwidth}{!}{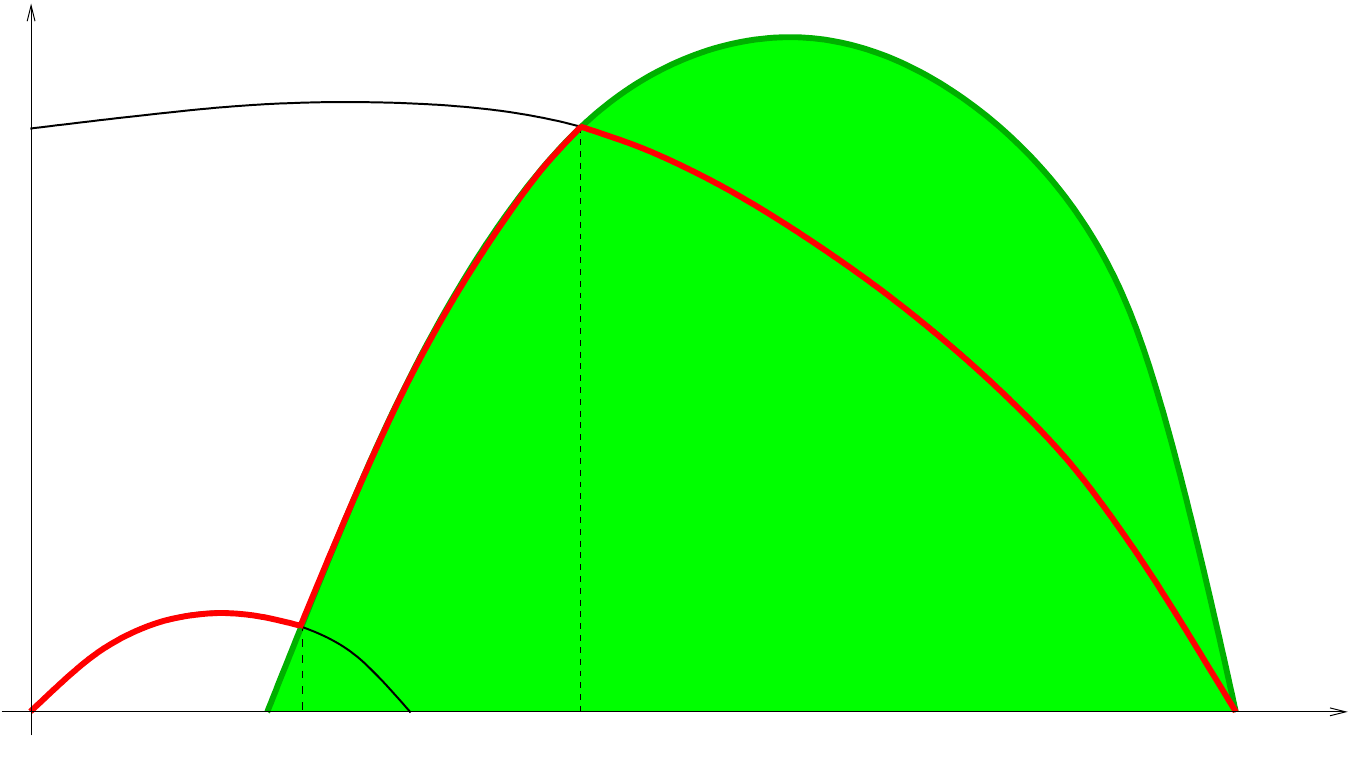} %\includegraphics[width=1\linewidth]{optprof.pdf}
\caption{The structure of the optimal solution under Condition (\ref{mon}).}
\label{Fig:optproof_1}
\end{figure}

Moreover, when $\beta^* \le \beta \le \beta^{**}$, the optimal profile $\gamma_\beta$ is the graph of the function (see Fig. \ref{Fig:optproof_1})
\begin{equation}
\label{optgraph}
P_\beta(U) =\begin{cases}
\Gamma_u(\beta, U) & U \in [0,U_u(\beta)], \\
P^*(U) & U \in [U_u(\beta), U_s(\beta)], \\
\Gamma_s(\beta, U) & U \in [U_s(\beta),1].
\end{cases}
\end{equation}
\begin{comment}
{\color{blue} where $U = U(x)$ is obtained by
\begin{equation*}
U_x = P_\beta(U)
\end{equation*}
and as in the previous section} $\Gamma_u(\beta, \cdot)$ and  $\Gamma_s(\beta, \cdot)$ {\color{blue} \sout{solve} are the unstable manifold of $(0,0)$ and the stable manifold of $(1,0)$, i.e. they solve} the uncontrolled ODE
\begin{equation}
\label{functeq}
P' + \frac{f(U)}{P} + \beta = 0
\end{equation}
with boundary conditions
\begin{equation*}
\Gamma_u(\beta, 0) = 0 \quad \text{and} \quad \Gamma_s(\beta, 1) = 0.
\end{equation*}
{\color{blue} These conditions corresponds to
\begin{equation*}
\lim_{x \to - \infty} U(x) = 0, \quad \lim_{x \to \infty U} U(x) = 1.
\end{equation*}
}
\end{comment}
In particular, there is an explicit expression for the optimal control, obtained by Equation \eqref{Equa:funct_eq_UP} for $P(U)$:
\begin{equation}
\label{Equa:alpha_beta}
\alpha_\beta(U)= \begin{cases}
\frac{F(U)+\beta}{U} & U \in [U_u(\beta),U_s(\beta)],\\
0 & U \in [0,U_u(\beta)]\cup[U_s(\beta), 1],
\end{cases}
\end{equation}
with
\begin{equation*}
%\label{F}
F(U) = (P^*)'(U) + \frac{f(U)}{P^*(U)} = \frac{3f(U)+Uf'(U)}{2\sqrt{Uf(U)}},  \quad \forall U \in (U^*,1).
\end{equation*}
This function also plays a fundamental role in the general case, i.e. when Assumption \eqref{mon} does not holds: indeed the level sets $\{F > - \beta\}$ describe the regions where the trajectories of the ODE \eqref{functeq} are entering the region $\{P < P^*\}$, and thus a control keeping $P(U) = P^*(U)$ can be applied.

%{\color{blue} 
Finally, even when $\beta>\beta^{**}$, a simple approximation argument \cite{bressancetraro} shows that the \emph{effort function $E(\beta)$} 
\begin{align*}
%\label{minfunct}
\begin{split}
E(\beta)=\inf \bigg\{ \int_0^1 \frac{1}{U} \left(P' +\frac{f(U)}{P} + \beta \right) \, dU&, P \text{ Lipschitz such that } \\
& \ P(0)=P(1)=0 \text{ and } P' +\frac{f(U)}{P} + \beta \ge 0 \bigg\},
\end{split}
\end{align*}
is equal to
\begin{equation*}
E(\beta) = \int_{U_u(\beta)}^{U_s(\beta)} \frac{F(U) + \beta}{U} dU.
\end{equation*}
We remark here that
\begin{equation*}
%\label{Equa:bound_integr}
F(U) \sim \begin{cases}
\frac{\sqrt{f'(U^*)}}{(U - U^*)^{\frac{1}{2}}} & U \searrow U^*, \\
- \frac{\sqrt{-f'(1)}}{(1 - U)^{\frac{1}{2}}} & U \nearrow 1,
\end{cases}
\end{equation*}
so that integrability is assured. In particular
\begin{equation*}
%\label{Equa:limi_F}
\lim_{U \to U^*} F(U) = + \infty, \quad \lim_{U \to 1} F(U) = - \infty.
\end{equation*}

The key role of Assumption (\ref{mon}) is that it implies a strict monotonicity condition for $F$ (i.e., it is strictly decreasing), which, in turn, guarantees that the upper level sets
$$
\{U\in [U^*,1]: \ F(U)+\beta>0\}
$$
are a single interval $(U^*,\check U)$. By observing Fig. \ref{Fig:optproof_1}, one sees that the intersection points of $\Gamma_u(\beta),\Gamma_s(\beta)$ with $P^*$ are unique, and that $\alpha_\beta$ as defined in \eqref{Equa:alpha_beta} is positive so that $\gamma_\beta$ is admissible.

\begin{remark}
\label{Rem:general_mon}
More generally, the same expression \eqref{optgraph} of $P_\beta$ is admissible if $F(U) + \beta$ has a single $0$: indeed, necessarily in the interval $(U_u(\beta),U_s(\beta))$ it holds $F(U) + \beta > 0$.
\end{remark}

\subsection{General case}
\label{Ss:general_case}

The analysis of the general case is performed in \cite{bianchini2026travelingprofilescontrolcost}. We first summarize the main results about uniqueness and structure of optimal traveling profiles.

\begin{definition}
\label{candidate}
Let $\beta \ge \beta^*$. We say that an a.c. profile $P: [0,1] \to [0,+\infty)$ is a \emph{candidate minimizer} for $\beta$, and write $P \in \tilde{\mathcal{A}}_\beta$, if the following holds:
\begin{enumerate}
\item $P(U) = \Gamma_u(\beta,U)$ in $[0,U_u(\beta)]$ and $P(U) = \Gamma_s(\beta,U)$ in $[U_s(\beta),1]$;
%\item if $\beta \geq \beta^{**}$, then $P(U) = P^*(U)$ in $[U^*,\gls{Ubar}]$;
\item $\Gamma_u(\beta,U) \leq P(U) \leq \Gamma_s(\beta,U)$;
%\item $P(U_u(\beta))=P^*(U_u(\beta)), \ P(U_s(\beta)=P^*(U_s(\beta))$;
\item for a.e. $U \in [U_u(\beta), U_s(\beta)]$
$$
\alpha(U) = \frac{1}{U} \left(P'+\frac{f(U)}{P}+\beta\right) \ge 0;
$$
\item $\alpha(U)=0$ for all $U \in \{U\in[U_u(\beta), U_s(\beta)]: P^*(U)\neq P(U)\}$.
\end{enumerate}
\end{definition}

Observe that the main point of the above definition is the requirement that the control is acting only on the critical curve $P^*(U) = \sqrt{U f(U)}$.

The following Structure Theorem characterizes the optimal traveling profiles and shows the key role of the set of candidates minimizers, introducing also the control set $\mathcal O(P_\beta)$ (the set where the control $\alpha$ is strictly positive) and the free set (where $P_\beta$ solves Equation \ref{functeq} with no control), and giving necessary and sufficient conditions for optimality in terms of a functional derivative equal to $0$.

\begin{theorem}
\label{Theo:regul_optP_intro}
Any optimal traveling profile $P_\beta$ is a candidate minimizer: the control $\alpha$ is applied only when $P_\beta(U) = P^*(U)$. In particular, $P_\beta$ is uniformly a.c. and Lipschitz if $U \not= U^*$.

Moreover, there exists an open set (the \emph{control set})
$$
\mathcal O(P_\beta) \subset \big\{ \Gamma_u < P^* < \Gamma_s \big\} \cap \mathrm{int} \{P_\beta = P^*\}
$$
where $\alpha > 0$. If $I = (U^-,U^+)$ is a connected component of the \emph{free set}
$$
\mathcal I(P_\beta) = \mathrm{int} \big( [0,1] \setminus \mathcal O(P_\beta) \big),
$$
then any optimal solution $P_\beta$ must satisfy the \emph{first-order conditions}
\begin{subequations}
\label{Equa:nec2_gen_pro_intro}
\begin{equation}
\label{Equa:nec2-_intro}
\text{if $P(U^-) > \Gamma_u(U^-)$} \quad \text{then} \quad \int_{U^-}^{\tilde U} \frac{Uf(U)-P^2(U)}{UP^2(U)} \, dU \ge 0 \quad \forall \tilde U \in I,
\end{equation}
%and
\begin{equation}
\label{Equa:nec2+_intro}
\text{if $P(U^+) < \Gamma_s(U^+)$} \quad \text{then} \quad \int_{\tilde U}^{U^+} \frac{Uf(U)-P^2(U)}{UP^2(U)} \, dU \le 0 \quad \forall \tilde U\in I.
\end{equation}
\end{subequations}

Finally, there exists a unique candidate minimizer $P_\beta$ satisfying the above first-order conditions, and the control set $\mathcal O(P_\beta)$ is monotone increasing w.r.t. $\beta$, and it consists of finitely many intervals whose boundary points depend smoothly on $(\beta,f)$ whenever $(\beta,f)$ belongs to on a open dense set of
$$
\mathcal D = \Big\{ (\beta,f), f \ \text{satisfies Assumption \eqref{H1}}, \ \beta > \beta^*(f) \Big\}.
$$ 
\end{theorem}

When studying the linearized ODE, it is useful to define the kernel
\begin{align}
\label{k}
K(\beta,U_1,U_2) = \int_{U_1}^{U_2} \frac{U f(U) - P_\beta(U)^2}{U P_\beta(U)^2} \, dU, \quad U_u(\beta) \leq U_1,U_2 \leq U_s(\beta). 
%\begin{cases}
%e^{\int_{U_1}^{U_2} \frac{U f(U) - P_\beta(f,U)^2}{U P_\beta(f,U)^2} \, dU} & \text{if} \ \beta^*(f) < \beta < \beta^{**}(f) \ \text{or} \ U^* \notin (U_1,U_2), \\
%0 & 0 < U_1 \leq \text{if} \ U^* < U_2 < 1 \ \text{and} \ \beta \geq \beta^{**}(f).
%\end{cases}
\end{align}
The conditions \eqref{Equa:nec2_gen_pro_intro} implies that if $I = (U^-,U^+)$ is a component of $\mathcal I(P_\beta)$ not containing $U_u(\beta),U_s(\beta)$, then
\begin{equation}
\label{Equa:kernel_0}
K(\beta,U^-,U^+) = 0.
\end{equation}

The next statement instead focuses on the regularity of the effort function $\beta,f \mapsto E(\beta,f)$. For our purposes, the key aspects are the Lipschitz dependence w.r.t. $\beta,f$ and the existence of an asymptote as $\beta \to \infty$.

\begin{theorem}
\label{Theo:Ebeta_intro}
The map $(\beta,f) \mapsto E(\beta,f)$ is Lipschitz in the domain of admissible sources and admissible speeds $\beta$
\begin{equation*}
\mathcal D = \Big\{ (\beta,f), f \ \text{satisfies Assumption \eqref{H1}}, \ \beta > \beta^*(f) \Big\} \subset \R \times C^2([0,1]).
\end{equation*}
Moreover
\begin{equation*}
\beta_1 \leq \beta_2, \ f_1 \leq f_2 \quad \Rightarrow \quad E(\beta_1,f_1) \leq E(\beta_2,f_2).
\end{equation*}
Moreover, the function $\beta \mapsto E(\beta)$ is $C^1$ in $\beta$, strictly increasing and
\begin{equation}
\label{Equa:asympt}
\lim_{\beta \to \infty} \bigg[ E(\beta) - \ln \bigg( \frac{1}{U^*} \bigg) \beta \bigg] = \int_{U^*}^1 \frac{F(U)}{U} dU = \int_{U^*}^1 \frac{2 \sqrt{f(U)}}{U^{3/2}} dU > 0.
\end{equation}
Finally, $f \mapsto E(\beta,f)$ is locally Lipschitz in the domain
\begin{equation*}
\mathcal D = \Big\{ (\beta,f), f \ \text{satisfies Assumption \eqref{H1}}, \ \beta > \beta^*(f) \Big\}.
\end{equation*} 
\end{theorem}

We finally recall the following explicit expression of the first derivative $\frac{dE}{d\beta}$ obtained in \cite{bianchini2026travelingprofilescontrolcost}: define the first-order perturbations of the optimal profile $P_\beta$
%\begin{subequations}
%\label{Equa:delta_P_intro}
\begin{align*}
%\label{Equa:delta_Pu_inrto}
\begin{cases}
\frac{d}{dU}\left( \frac{\partial P_u(U)}{U} \right) = \frac{Uf(U)-P(U)^2}{UP(U)^2} \left(\frac{\partial P_u(U)}{U} \right) - \frac{1}{U}, \\ %- \frac{\delta f(U)}{P(U)}, \\
\partial P_u(U_u(\beta)) = \partial_\beta \Gamma_u(\beta,U_u(\beta)),
\end{cases}
\end{align*}
\begin{align*}
%\label{Equa:delta_Ps_inrto}
\begin{cases}
\frac{d}{dU}\left( \frac{\partial P_s(U)}{U} \right) = \frac{Uf(U)-P(U)^2}{UP(U)^2} \left(\frac{\partial P_s(U)}{U} \right) - \frac{1}{U}, \\ %- \frac{\delta f(U)}{P(U)},\\
\partial P_s(U_s(\beta)) = \partial_\beta \Gamma_s(\beta,U_s(\beta)).
\end{cases}
\end{align*}
%\end{subequations}
The function $\partial P_u$ is the perturbation of the solution by varying $\beta$ but keeping the control the same, with initial condition $\partial P_u(U = 0) = 0$; the same holds for $\partial P_s(U)$, with the end condition $\partial P_s(U=1) = 0$. The existence of these functions and their properties are proved by tools of Dynamical System Theory, we refer to \cite{bianchini2026travelingprofilescontrolcost}. It can be proved that the explicit expression of the derivative is
\begin{equation*}
%\label{Equa:explit_deri}
\frac{d}{d\beta} E(\beta) = %\begin{cases}
\frac{\partial P_s(U,\beta) - \partial P_u(U, \beta)}{U} \quad \forall U \in \clos(\mathcal O(P_\beta)). \\
%\frac{\partial P_s(U,\beta) - \partial P_u(U,\beta)}{U} e^{K(\beta,U,U^+)} & U \in (U^-,U^+) \in \gls{IcalPopt} \ \text{connected component}.
%\end{cases} %[U_u(\beta),U_s(\beta)]. %\overline{\mathcal{O}(\beta)}.
\end{equation*}

Similar computations can be performed when dealing with the second derivative, under the assumption that the free set $\mathcal I(P_\beta)$ is made by finitely many intervals $I_i = (U^-_i(\beta),U^+_i(\beta))$, $i=1,\dots,N$. Let
\begin{equation*}
\hat U_u(\beta) = \max \big\{ U : \Gamma_u(\beta) = P_\beta(U) \big\}, \quad \hat U_u(\beta) = \min \big\{ U : \Gamma_u(\beta) = P_\beta(U) \big\}.
\end{equation*}
The explicit expression of the second derivative is then
\begin{equation}
\label{Equa:second_der_E_beta}
\begin{split}
\partial^2_\beta E(\beta) &= \frac{1}{\hat U_s(\beta)} \frac{\partial^2 \Gamma_s}{\partial \beta^2}(\beta,\hat U_s(\beta)) - \frac{1}{\hat U_u(\beta)}\frac{\partial ^2 \Gamma_u}{\partial \beta^2}(\beta,\hat U_u(\beta)) \\
& \quad + \sum_{i=1}^N \int_{U_i^-(\beta)}^{U_i^+(\beta)} \frac{2f(W)}{WP_i^3(W, \beta)} \left(\frac{\partial P_i}{\partial \beta}(W, \beta)\right)^2\, e^{-K(\beta, U_i^-(\beta), W)} dW.
\end{split}
\end{equation}
where $K(\beta,U_1,U_2)$ is the kernel defined in \eqref{k} and $P_i = P_i(U,\beta)$ are the solution to the ODE \eqref{functeq} in the $i$-th connected component $I_i(\beta)$:
\begin{equation*}
%\label{functeq}
\frac{dP_i}{dU} = - \frac{f(U)}{P_i(U)} - \beta, \quad P_i(U^-_i(\beta)) = P^*(U^-_i(\beta)) \ (\text{or equivalently} \ P_i(U^+_i(\beta)) = P^*(U^+_i(\beta))).
\end{equation*}
As we stated above, a key part of the computation of the derivative is the fact that $U^\pm_i(\beta)$ depends smoothly on $\beta$: this is not the case when a splitting occurs in one of the intervals. We will show that, in this situation, the second derivative generally does not exist.

Before ending, we want to observe that in the case $P_\beta$ has the structure of \eqref{optgraph}, the formula for the second derivative becomes simply
\begin{equation}
\label{Equa:second_simply}
\partial^2_\beta E(\beta) = \frac{1}{\hat U_s(\beta)} \frac{\partial^2 \Gamma_s}{\partial \beta^2}(\beta,\hat U_s(\beta)) - \frac{1}{\hat U_u(\beta)}\frac{\partial ^2 \Gamma_u}{\partial \beta^2}(\beta,\hat U_u(\beta)).
\end{equation}

\section{Monotonicity condition}
\label{S:monoton_cond}

\begin{figure}
\begin{subfigure}{.475\textwidth}
\resizebox{\textwidth}{!}{\includegraphics{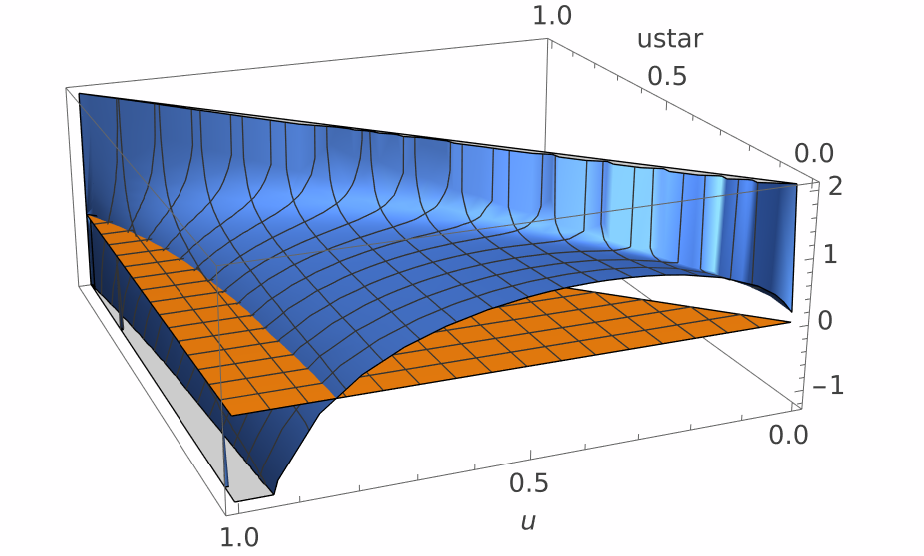}}
\end{subfigure} \hfill
\begin{subfigure}{.475\textwidth}
\resizebox{\textwidth}{!}{\includegraphics{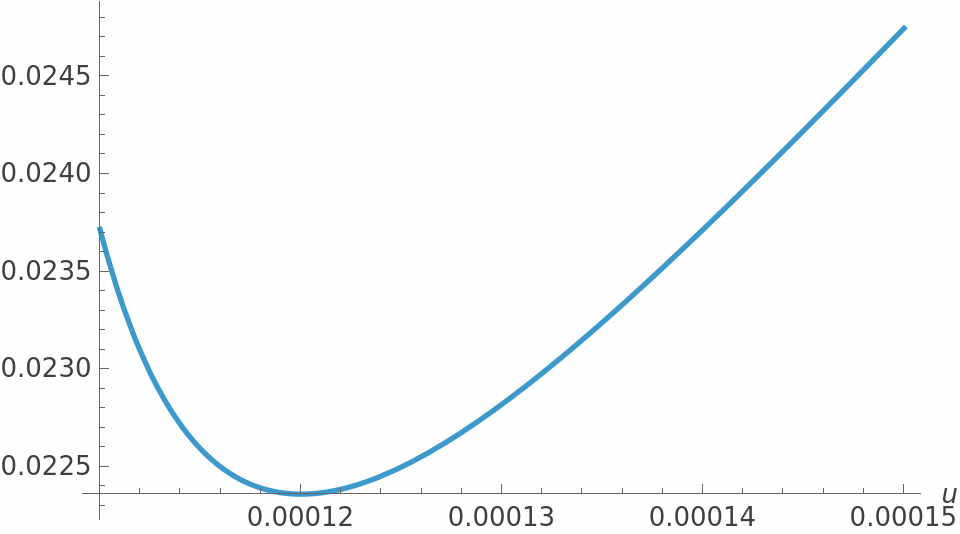}}
\end{subfigure}
\caption{The function $F(U)$ for the cubic \eqref{Equa:f_cubic} in the region $U^* < U < 1$. %Observe that $F$ is not monotone for $U^*$ close to $0$, even if \eqref{adm} may still hold. 
On the right the plot of $F(U)$ for $U^* = 10^{-4}$: the minimum can be shown to be at $\sim \frac{6 U^*}{5}$. %, so that the manifold $\Gamma_u$ intersects $P^*$ at $\sim 10^{-5}$: this point is on the right of the minimum $1.058 * 10^{-5}$.
}
\label{Fig:fcubicF}
\end{figure}

In this section we numerically compute the function $F(U)$ for the cubic source
\begin{equation}
\label{Equa:f_cubic}
f(U) = U (U - U^*) (1 - U).
\end{equation}
The function $F$ is represented in Fig. \ref{Fig:fcubicF} (left) as a function of $U^*,U$. Note that the function is not monotone for $U^* < 0.4$: for $U^* \ll 1$ there are two critical points, a local minimum and a local maximum. In particular, for $0 < -\beta \ll 1$ there are $3$ solution of $F(U) + \beta = 0$:
\begin{equation*}
\{F(U) + \beta > 0\} = (U^*,U_1) \cup (U_2,U_3), \quad U_1 < U_2,\ U_3 < 1.
\end{equation*}
We thus conclude that the Condition (\ref{mon}) fails even for the cubic source. 

The value $U^* = 0.4$ can be computed explicitly: indeed, it corresponds to the value such that there exists $U \in (0,1)$ s.t.
\begin{equation*}
\frac{d}{dU} F(U) = 0, \quad \frac{d^2}{dU^2} F(U) = 0.
\end{equation*}
For a deeper analysis of this case we refer to \cite{bianchini2026travelingprofilescontrolcost}, where more general polynomial sources are considered.

\section{$E(\beta)$ not convex}
\label{S:E_not_convex}

As a first step, we compute the explicit formulas for the stable and unstable manifolds $\Gamma_u,\Gamma_s$ and their derivatives w.r.t. $\beta$ in the case $\beta = 0$: this case is easier because the ODE is now Hamiltonian. Since in our construction the monotonicity assumption holds , the effort is computed explicitly by
\begin{equation*}
E(\beta) = \int_{U_u(\beta)}^{U_s(\beta)} \frac{F(U) + \beta}{U} dU,
\end{equation*}
and its second derivative is given by \eqref{Equa:second_simply}.

\subsubsection{Formulas for $\Gamma_u,\Gamma_s$ in the case $\beta = 0$}
\label{Ss:formulas_beta0}

Fix $\beta=0$ and define the primitive function
\begin{equation}
\label{Equa:G_primi_f}
G(U) = - \int_0^U f(V) dV.
\end{equation}
Note that, clearly, $G(U)> 0$ for all $U \in [0, U^*]$. Moreover, by assumption
\begin{align*}
G(1) = - \int_0^1 f(V) dV < 0,
\end{align*}
so that we can find $\bar U \in (U^*,1)$ such that $G(\bar U)=0$ by the intermediate value theorem. Finally, since
\begin{equation*}
G'=-f \le 0 \quad \text{in} \ [U^*,1],
\end{equation*}
we also deduce
\begin{align*}
G(U)\ge G(1) \quad \forall U \in [U^*,1].
\end{align*}
In this case, Equation \eqref{functeq} is an Hamiltonian system and can be rewritten as
\begin{equation*}
\frac{d}{dU} \left(\frac{P^2(U)}{2}\right) = -f(U) = G'(U), 
\end{equation*}
which implies, by the previous observations,
\begin{equation}
\label{unstable0}
\Gamma_u(0,U) = \sqrt{2G(U)} \quad \forall U \in [0,\bar U],
\end{equation}
%and
\begin{equation}
\label{stable0}
\Gamma_s(0,U) = \sqrt{2 \left( G(U)-G(1) \right)} \quad \forall U\in [0,1].
\end{equation}
Moreover, for $\beta=0$, %the equation for the derivative $\partial_\beta \Gamma_u$ becomes %\ref{fder} becomes
\begin{align}
\label{delta}
\frac{\partial}{\partial U} \left( \Gamma_u(\beta, U) \frac{\partial \Gamma_u}{\partial \beta}(\beta, U) \right) = - \Gamma_u(\beta,U),
\end{align}
and the same for $\Gamma_s$. By the boundary conditions
\begin{align*}
\Gamma_u(0,0) = \Gamma_s(1,0)=0,
\end{align*}
we obtain the explicit formulas
\begin{equation}
\label{eqmeno}
\frac{\partial \Gamma_u}{\partial \beta}(0, U)=-\frac{1}{\Gamma_u(0, U)}\int_0^U \Gamma_u(0, U') \, dU'= -\frac{1}{\sqrt{G(U)} }\int_0^U \sqrt{G(U')} dU',
\end{equation}
%    and
\begin{equation*}
%\label{eqplus}
\frac{\partial \Gamma_s }{\partial \beta}(0, U) = -\frac{1}{\Gamma_s(0,U)} \int_1^U \Gamma_s(0, U') \, dU'= \frac{1}{\sqrt{G(U) - G(1)}} \int_U^1 \sqrt{G(U') - G(1)} dU'.
\end{equation*}
Similarly, the equation for $\partial^2_\beta \Gamma_u$ is for $\beta = 0$ %by \ref{sder} evaluated at $\beta=0$,
\begin{equation}
\label{deltas}
\frac{\partial}{\partial U}
\left(\Gamma_u(0, U) \frac{\partial ^2 \Gamma_u }{\partial \beta^2} (0, U) \right) = - 2 \left[ \frac{\partial}{\partial U} \left(\frac{\partial \Gamma_u}{\partial \beta} \right)(0, U) + 1 \right] \frac{\partial \Gamma_u}{\partial \beta}(0, U),
\end{equation}
and the same for $\Gamma_s$. Thus,
\begin{equation}
\label{comp}
\frac{\partial ^2 \Gamma_u }{\partial \beta^2} (0, U)=\frac{1}{\Gamma_u(0, U)}   \left[-\left(\frac{\partial \Gamma_u}{\partial\beta}(0,U)\right)^2 -2 \int_0^U  \frac{\partial \Gamma_u}{\partial \beta}(U',0) \, dU' \right],
\end{equation}
%and
\begin{equation*}
%\label{comps}
\frac{\partial ^2\Gamma_s}{\partial \beta^2}(U,0)=\frac{1}{\Gamma_s(0,U)}    \left[-\left(\frac{\partial \Gamma_s}{\partial\beta}(U,0)\right)^2 +2 \int_U^1  \frac{\partial \Gamma_s}{\partial \beta}(U',0) \, dU' \right].
\end{equation*}
We remark that for all the above computation to be meaningful, we only need $f$ to be continuous: by approximation, the second derivatives $\frac{\partial ^2 \Gamma_u}{\partial \beta^2}$, $\frac{\partial ^2 \Gamma_s}{\partial \beta^2}$ remain close to the value we compute below, and then the sign of \eqref{Equa:second_simply} will still be negative. %the same non-convexity result holds for any smooth source sufficiently close in $C^0$ to the $f$ we construct below (recall also the Lipschitz regularity of $E(\beta,f)$ stated in Theorem \ref{Theo:Ebeta_intro}).

\subsubsection{Definition of the source $f$}
\label{Sss:explicitc_f}

Let $K \ge 1$ and $U^* \in (0, \frac{1}{2})$ and define the function
\begin{align}\label{Eq:nonconvsource}
f(U)=f_K(U) = \begin{cases}
U(U-U*) & 0\le U \le U^*, \\
K(U-U^*) & U^* \le U \le  \frac{1}{2}, \\
K(1-2U^*)(1-U) & \frac{1}{2}\le U \le 1.
\end{cases}
\end{align}
We compute
\begin{equation*}
G(U) = \begin{cases}
U^*\frac{U^2}{2}-\frac{U^3}{3} & 0\le U \le U^* \\
G(U^*)-\frac{K}{2}\left(U-U^*\right)^2 & U^* \le U \le  \frac{1}{2} \\
G(\frac{1}{2}) +\frac{K}{2}(1-2U^*) [ (U-1)^2 - \frac{1}{4} ] & \frac{1}{2}\le U \le 1,
\end{cases}
\end{equation*}
with $G(U^*)=\frac{1}{6}(U^*)^3$ and $G(\frac{1}{2})=G(U^*) -\frac{K}{2}(\frac{1}{2}-U^*)^2.$ Note, in particular, that condition $U^* <\frac{1}{2}$ ensures
\begin{equation*}
G(1) =-\int_0^1 f(U) \,dU =\frac{1}{6}\left(U^*\right)^3-\frac{K}{2} (1-U^*) \bigg( \frac{1}{2} - U^* \bigg) < 0 \quad \text{ for } K \gg 1,
\end{equation*}
that is $f$ is positively predominant.

It is fairly easy to see that $F(U)$ has a single $0$, which gives that for $\beta$ close to $0$ the structure of the optimal profile $P_\beta$ is as in \eqref{optgraph} (see Remark \ref{Rem:general_mon}): using \eqref{Eq:nonconvsource} and the formula for $F(U)$ one verifies that $F(U) = 0$ only for $U = \frac{3}{4}$, and moreover $\frac{d}{dU} F(U=\frac{3}{4}) = -4$, so that this condition is stable for speeds $\beta$ close to $0$.

\subsubsection{Estimates on the manifold $\Gamma_u$}
\label{Sss:couter_esimt_Gammau}

For $0\le U\le U^*$, by (\ref{eqmeno}), we can directly compute
\begin{align*}
\frac{\partial\Gamma_u}{\partial\beta}(0,U) &= -\frac{1}{U\sqrt{1-\frac{2U}{3U^*}}} \int_0^U U' \sqrt{1-\frac{2U'}{3U^*}} \, dU' \\
&= - \frac{(\frac{3U^*}{2})^2}{U \sqrt{1 - \frac{2U}{3U^*}}} \bigg[ \frac{2}{3} \bigg( 1 - \bigg( 1 - \frac{2U}{3U^*} \bigg)^\frac{3}{2} \bigg) - \frac{2}{5} \bigg( 1 - \bigg( 1 - \frac{2U}{3U^*} \bigg)^\frac{5}{2} \bigg) \bigg] \\
&= - \frac{\frac{2 U^*}{5}}{\sqrt{1-\frac{2U}{3U^*}} ( 1 + \sqrt{1-\frac{2U}{3U^*}} )} + \frac{U^*}{5} - \frac{2U}{5}.
\end{align*}
In particular,
\begin{equation*}
\frac{\partial \Gamma_u}{\partial\beta}(0,U^*) = - \bigg( \frac{3\sqrt{3}-2}{5} \bigg) U^*.
\end{equation*}
It follows that, for $0\le U \le U^*$,
\begin{align*}
\int_0^{U} \frac{\partial\Gamma_u}{\partial\beta}(0,U') \, dU' &= - \frac{2U^*}{5} \int_0^U (-3U^*) \frac{d}{dU'} \bigg[ \ln \bigg( 1 + \sqrt{1 - \frac{2U'}{3U^*}} \bigg) \bigg] \, dU' + \frac{U^*U}{5} - \frac{U^2}{5} \\
&= - \frac{6 (U^*)^2}{5} \ln \bigg( \frac{2}{1 + \sqrt{1 - \frac{2U}{3U^*}}} \bigg) + \frac{U^* U}{5} - \frac{U^2}{5},
\end{align*}
and, consequently,
\begin{equation*}
\int_0^{U^*} \frac{\partial \Gamma_u}{\partial\beta}(0, U') \, dU' = - \frac{6 (U^*)^2}{5} \log{(3-\sqrt{3})}.
\end{equation*}
By identity (\ref{comp}), we conclude
\begin{align*}
\frac{\partial ^2 \Gamma_u}{\partial \beta^2}(0,U^*) &= \frac{1}{\sqrt{2G(U^*)}} \bigg( - \bigg( \frac{3\sqrt{3} - 2}{5} \bigg)^2 (U^*)^2 + \frac{12}{5} \ln(3-\sqrt{3}) (U^*)^2 \bigg) \\
&= \frac{\sqrt{3U^*}}{5} \bigg( 12 \ln(3-\sqrt{3}) - \frac{31-12\sqrt{3}}{5} \bigg) > 0.
\end{align*}

Thanks to formulas (\ref{unstable0}) and (\ref{stable0}), we can explicitly compute the unique intersections with the critical curve $P^*$
\begin{equation*}
U_u(0) = \frac{3}{4} U^* + \sqrt{ \bigg( \frac{U^*}{4} \bigg)^2 + \frac{G(U^*)}{K}} \quad \text{and} \quad U_s(0) = \frac{1}{2}.
\end{equation*}
Note, in particular, that, by the inequality
\begin{equation*}
\sqrt{1+x}-1 \le \frac{x}{2} \quad \forall x\ge 0,
\end{equation*}
we obtain
\begin{equation}
\label{dist}
0 \le U_u(0) - U^* = - \frac{U^*}{4} + \sqrt{\bigg( \frac{U^*}{4} \bigg)^2 + \frac{G(U^*)}{K}} \le \frac{(U^*)^2}{3K}.
\end{equation}

For every $U \in [U^*, \bar U]$ we can write
\begin{equation*}
\Gamma_u(0,U) = \Gamma_u(0, U^*) - R_1(U),
\end{equation*}
with
\begin{equation*}
0 \le R_1(U)= \sqrt{2G(U^*)}- \sqrt{2G(U^*) -K(U-U^*)^2} \le \sqrt{K}(U-U^*).
\end{equation*}
In particular
\begin{equation*}
R_1(U_u) \leq \frac{(U^*)^2}{3 \sqrt{K}}.
\end{equation*}

Moreover, by \eqref{delta}, there holds
\begin{equation*}
\frac{\partial \Gamma_u}{\partial \beta}(0,U)=\frac{\partial \Gamma_u}{\partial \beta}(0,U^*)+ R_2(U)
\end{equation*}
with
\begin{align*}
R_2(U) = - (U - U^*) + \frac{1}{\Gamma_u(0, U^*)} \int_{U^*}^U R_1(U') \, dU' + \frac{\partial \Gamma_u}{\partial \beta}(0,U) \frac{R_1(U)}{\Gamma_u(0, U^*)},
\end{align*}
and, thus, for all $U^* \le U \le U_u(0) \le \bar U$ and $K$ sufficiently large, by \eqref{dist},
\begin{align*}
\left|R_2(U)\right|&\le (U-U^*)+ \frac{1}{2\Gamma_u(0, U^*)}\sqrt{K} (U-U^*)^2 + \frac{C}{\Gamma_u(0, U^*)}\sqrt{K}(U-U^*) \le \frac{C}{\sqrt{K}},
\end{align*}
For simplicity, throughout the remainder of this subsection, $C$ is a numeric constant whose value only depends on $U^*$ and may vary from time to time

Similarly, by (\ref{deltas}), we can write
\begin{align*}
\Gamma_u(0,U) \frac{\partial ^2 \Gamma_u}{\partial \beta^2}(0, U)=\Gamma_u(0,U^*)\frac{\partial ^2 \Gamma_u}{\partial \beta^2}(0,U^*)+ \frac{\partial \Gamma_u}{\partial \beta}(0, U^*)^2 -\frac{\partial \Gamma_u}{\partial \beta}(0,U)^2 -2\int_{U^*}^U \frac{\partial \Gamma_u}{\partial \beta} (0, U') \, dU',
\end{align*}
which implies
\begin{equation*}
\frac{\partial ^2 \Gamma_u}{\partial \beta^2}(0,U)= \frac{\partial ^2 \Gamma_u }{\partial \beta^2}(0, U^*)+R_3(U),
\end{equation*}
with
\begin{align*}
R_3(U) :&= \frac{1}{\Gamma_u(0,U^*)}\left[\left( \frac{\partial \Gamma_u}{\partial \beta}(0,U^*)\right)^2 -\left(\frac{\partial \Gamma_u}{\partial \beta}(0,U)\right)^2\right] \\
& \quad -\frac{2}{\Gamma_u(0, U)}\int_{U^*}^U \frac{\partial \Gamma_u}{\partial \beta}(0, U') \, dU'\\
& \quad - \frac{\partial ^2\Gamma_u}{\partial \beta^2}(0,U^*)\left(\frac{\Gamma_u(0, U^*)}{\Gamma_u(0,U)}-1 \right)
\end{align*}
and, therefore, also
\begin{equation*}
\left|R_3(U)\right| \le \frac{C}{\sqrt{K}}, \quad \forall U^*\le U \le U_u(0) \text{ and } K \gg 1.
\end{equation*}

In particular, we conclude that, for $K$ sufficiently large,
\begin{equation*}
\frac{\partial ^2 \Gamma_u }{\partial \beta^2}(0,U_u(0)) \ge \frac{1}{2}\frac{\partial ^2 \Gamma_u}{\partial \beta^2}(0,U^*) > 0.
\end{equation*}

\subsubsection{Estimates on the manifold $\Gamma_s$}
\label{Sss:estim_couteGammas}

Conversely, equations (\ref{delta}) and (\ref{deltas}) for $\Gamma_s$ can be directly integrated to find
\begin{equation*}
\Gamma_s(0,U) = \sqrt{K(1-2U^*)}(1-U), \quad \frac{\partial \Gamma_s}{\partial \beta}(0,U) = \frac{1-U}{2}, \quad \frac{\partial ^2\Gamma_s}{\partial \beta^2}(0,U)= \frac{1-U}{4\sqrt{K(1-2U^*)}}.
\end{equation*}

\subsubsection{Computing the second derivative of $E(\beta)$ for $\beta = 0$ and non-convexity of $E$}
\label{Sss:comput_second}

In particular,
\begin{equation*}
\frac{\partial^2\Gamma_s}{\partial \beta^2}(0, U_s(0))=\frac{1}{8\sqrt{K(1-2U^*)}} \to 0 \text{ as } K\to +\infty.
\end{equation*}
We therefore conclude that it is possible to pick $K$ sufficiently large so that
\begin{equation*}
\partial^2_\beta E(0) = \frac{1}{U_s(0)} \frac{\partial^2 \Gamma_s}{\partial \beta^2}(0,U_s(0)) - \frac{1}{U_u(0)}\frac{\partial^2 \Gamma_u}{\partial \beta^2}(0,U_u(0))<0.
\end{equation*}

%\begin{remark}
%\label{Rem:general}
%Since $E$ depends continuously from the source $f(U)$ as proved in Section \ref{S:regula_E_beta_f}, we conclude that the set of sources $f(U)$ such that $E$ is not convex is a nonempty open set in $C^2$.
%\end{remark}

\section{Non-convexity and anisotropic perimeter}
\label{S:anisotr_nonconvex}

In this section the non convexity is proved by checking that it is incompatible with the linear growth at $U \to +\infty$. We recall that by \eqref{Equa:asympt} there exists an asymptote as $U \to +\infty$.

In order to have that $E(\beta)$ corresponds to the integrand of an anisotropic perimeter functional, the following assumptions are made on the effort function $E(\beta)$:
\begin{enumerate}
\item $\beta \mapsto E(\beta)$ is convex,
\item $E(\beta) - \beta E'(\beta) \geq 0 \quad \forall \beta\ge0.$
\end{enumerate}
Because of the existence of the asymptote by Theorem \ref{Theo:Ebeta_intro}, the first assumption implies the second. Indeed,  if we assume convexity a priori, in particular $\beta \mapsto E'(\beta)$ is increasing, so that the asymptotic behavior of $E$ implies
\begin{equation*}
E'(\beta)\le\ln \bigg( \frac{1}{U^*} \bigg) \quad \forall \beta \in \mathbb R.
\end{equation*}
In particular, the function $\beta \mapsto E(\beta)-\beta \ln \big( \frac{1}{U^*} \big)$ is decreasing, so that, for all $\beta \ge0$,
\begin{equation}
\label{Eq:aboveas}
E(\beta)-\beta E'(\beta) \ge E(\beta) -\ln \bigg( \frac{1}{U^*} \bigg) \beta \ge \lim_{\beta \to \infty} \bigg[ E(\beta) - \ln \bigg( \frac{1}{U^*} \bigg) \beta \bigg] = \int_{U^*}^1 \frac{2 \sqrt{f(U)}}{U^{3/2}} dU >0.
\end{equation}
Conversely, similar computations show here that for the cubic source $f(U) = U (U - U^*) (1 - U)$ there are $U^*$'s such that it holds
\begin{equation*}
E(\beta = 0) < \int_{U^*}^1 \frac{2 \sqrt{f(U)}}{U^{\frac{3}{2}}} dU.
\end{equation*}
Hence $E(\beta)$ is not convex even for the cubic source.

For this purpose, just consider $U^* = \frac{1}{2}$, so that $E(\beta = 0) = 0$, since the system \eqref{sys} is Hamiltonian and $G(0) = G(1)$, where $G$ is the potential defined by \eqref{Equa:G_primi_f}. Observe that
\begin{equation*}
\int_{\frac{1}{2}}^1 \frac{2 \sqrt{f(U)}}{U^{\frac{3}{2}}} dU > 0.
\end{equation*}
Therefore, by the Lipschitz dependence of Theorem \ref{Theo:Ebeta_intro}, for $U^*$ close to $\frac{1}{2}$ and $\beta$ close to $0$, the effort $E(\beta)$ satisfies
\begin{equation*}
E(\beta) < \ln \bigg( \frac{1}{U^*} \bigg) \beta + \int_{U^*}^1 \frac{2 \sqrt{f(U)}}{U^{\frac{3}{2}}} dU,
\end{equation*}
which implies non-convexity (by contradicting \eqref{Eq:aboveas}).

\section{Example of $E$ with only $C^{1,1}$-regularity}
\label{S:exampl_C11}

In this section we give a mechanism forbidding the effort $E(\beta)$ to be of class $C^2$. A more refined analysis of trajectories near the point $U^*$ yields that $E(\beta)$ is not even $C^{1,1}$, see \cite{bianchini2026travelingprofilescontrolcost}. The phenomenon we want to analyze here is what happens when the topological structure of the control set $\mathcal O(P_\beta)$ varies (or equivalently the structure of the free set $\mathcal I(P_\beta)$ changes). In the same paper \cite{bianchini2026travelingprofilescontrolcost} it is shown that it is possible to construct a source $f(U)$ such that the optimal profile has an arbitrarily large number of oscillations: preliminary numerical analysis gives that the splitting described here occurs.

By Theorem \ref{Theo:regul_optP_intro}, the free set is made by countably many intervals $I(\beta) = (U^-(\beta),U^+(\beta))$, determined uniquely by the conditions \eqref{Equa:nec2-_intro}, \eqref{Equa:nec2+_intro}. The fact that the control set $\mathcal O(P_\beta)$ is increasing (in the sense of set inclusions) implies that the connected components $I(\beta)$ shrink or even split as $\beta$ increases. We want to prove here that the dependence of $E(\beta)$ is only $C^{1,1}$ in the case the interval $I(\beta)$ splits into two.

We consider $f$ is an admissible sources and we assume for simplicity that for $\beta < \bar \beta$
\begin{equation*}
I(\beta) = (U_1^-(\beta),U_1^+(\beta)),
\end{equation*}
whereas for $\beta > \bar{\beta}$
\begin{equation*}
I(\beta)=(U_2^-(\beta), U_2^+(\beta)) \cup (U_3^-(\beta), U_3^+(\beta)).
\end{equation*}
This situation is in accordance with the fact that the set $\mathcal I(P_\beta)$ is decreasing. We will compute the right and left limits of $\partial^2_\beta E$ at $\bar \beta$ by means of formula \eqref{Equa:second_der_E_beta}: more precisely we will compute the contribution of the term
\begin{equation*}
\int_{U_i^-(\beta)}^{U_i^+(\beta)} \frac{2f(W)}{WP_i^3(W, \beta)} \left(\frac{\partial P_i}{\partial \beta}(W, \beta)\right)^2\, e^{-K(\beta, U_i^-(\beta), W)} dW
\end{equation*}
to the second derivative, and see that it has a discontinuity.

For $i=1,2,3$, denote by $P_i=P_i(U, \beta)$ the solution to
\begin{equation*}
\begin{cases}
P_i'(\beta) + \frac{f(U)}{P_i} + \beta=0, \quad U \in (U_i^-(\beta), U_i^+(\beta)), \\
P_i(U_i^-(\beta)) = P^*(U_i^-(\beta)).
\end{cases}
\end{equation*}
For $\beta < \bar \beta$, the equation of the perturbation is
\begin{equation*}
%\label{Equa:deltaP_1_eq}
\frac{d}{dU} \delta P_1 = \frac{f(U)}{P^2} \delta P_1 - 1.
\end{equation*}
The initial condition $\delta P_1$ at some point $U \in (U^-_1,U^+_1)$ is computed in order to maintain
\begin{equation*}
\int_{U^-_1}^{U^+_1} \frac{U f(U) - P^2}{U P^2} dU = 0,
\end{equation*}
which becomes
\begin{equation}
\label{Equa:require_cc}
\int_{U^-_1}^{U^+_1} \frac{2 f(U) \delta P_1}{P^3} dU = 0.
\end{equation}
Note that this requires
\begin{equation}
\label{Equa:sing_deltaP}
\delta P_1(U_1^+) < 0 \quad \text{and} \quad \delta P_1(U^-_1) > 0,
\end{equation}
otherwise the solution would be strictly positive or negative in the interval, respectively.

For $\beta > \bar \beta$ the interval $(U^-_1,U^+_1)$ splits into two intervals and then the initial data for them is similarly computed: if $\delta P_2$, $\delta P_3$ are the solutions in the intervals $(U^-_2,U^+_2),(U^-_3,U^+_3)$, then
\begin{equation*}
\int_{U^-_2}^{U^+_2} \frac{2 f(U) \delta P_2}{P^3} dU = \int_{U^-_3}^{U^+_3} \frac{2 f(U) \delta P_3}{P^3} dU = 0.
\end{equation*}
In particular by \eqref{Equa:sing_deltaP} we have
\begin{equation}
\label{Equa:sign_deltaP_23}
\delta P_2(U_2^+) < 0 \quad \text{and} \quad \delta P_3(U^-_3) > 0.
\end{equation}

For $\beta = \bar \beta$, the optimal profile $P(U)$ is the same for both cases, the difference being in the initial data: in particular, we can write $\delta P_1$ in terms of $\delta P_2, \delta P_3$ and the solution of the homogeneous equation as
\begin{equation*}
%\label{Equa:delta_P_1_gener_23}
\delta P_1(\bar \beta,U) = \begin{cases}
a e^{\int_{U_m}^U \frac{f(W)}{P(W)^2} dW} + \delta P_2(\bar \beta,U) & U \in (U^-_1,U_m], \\
b e^{\int_{U_m}^U \frac{f(W)}{P(W)^2} dW} + \delta P_3(\bar \beta,U) & U \in (U_m,U_1^+),
\end{cases}
\end{equation*}
with
\begin{equation*}
%\label{Equa:ab_rel0}
a + \delta P_2(U_m) = b + \delta P_3(U_m),
\end{equation*}
Imposing \eqref{Equa:require_cc} we obtain also
\begin{equation*}
%\label{Equa:ab_rel}
a \int_{U_1^-}^{U_m} \frac{2 f(W)}{P^3} e^{\int_{U_m}^W \frac{f(V)}{P^2(V)} dV} dW + b \int_{U_m}^{U^+_1} \frac{2 f(W)}{P^3} e^{\int_{U_m}^W \frac{f(V)}{P^2(V)} dV} dW = 0.
\end{equation*}
We observe also that $\delta P_2(U_m) < 0$ and $\delta P_3(U_m) > 0$ by \eqref{Equa:sign_deltaP_23}, so that being the linear system above invertible we conclude that
\begin{equation*}
a = \frac{\int_{U_m}^{U^+_1} \frac{2 f(W)}{P^3} e^{\int_{U_m}^W \frac{f}{P^2} dV} dW}{\int_{U^-_1}^{U^+_1} \frac{2 f(W)}{P^3} e^{\int_{U_m}^W \frac{f}{P^2} dV} dW} \big( \delta P_3(U_m) - \delta P_2(U_m) \big) > 0,
\end{equation*}
\begin{equation*}
b = - \frac{\int_{U_1^-}^{U_m} \frac{2 f(W)}{P^3} e^{\int_{U_m}^W \frac{f}{P^2} dV} dW}{\int_{U^-_1}^{U^+_1} \frac{2 f(W)}{P^3} e^{\int_{U_m}^W \frac{f}{P^2} dV} dW} \big( \delta P_3(U_m) - \delta P_2(U_m) \big) < 0.
\end{equation*}

We now compute the second derivative of $E$ at $\bar \beta-$ using the above expression:
\begin{equation*}
\begin{split}
\partial^2_\beta E(\bar \beta) &= \int_{U_1^-(\bar \beta)}^{U_1^+(\bar \beta)} \frac{2f(W)}{W P(\bar \beta,W)^3} \bigg( \frac{\partial P_1}{\partial \beta}(W,\bar \beta) \bigg)^2 e^{-K(\bar \beta, U_1^-(\beta), W)} dW \\
&= \int_{U^-_1}^{U^+_1} \frac{2 f(W) e^{-K(\bar \beta, U_1^-, W)}}{W P(\bar \beta,W)^3} \Big( \big( a e^{\int_{U_m}^U \frac{f(W)}{W^2} dW} + \delta P_2 \big) \rest_{(U^-_1,U_m)} + \big( b e^{\int_{U_m}^U \frac{f(W)}{W^2} dW} + \delta P_3 \big) \rest_{(U_m,U^+_1)} \Big)^2 dW \\
&= \int_{U^-_1}^{U_m} \frac{2 f(W) e^{-K(\bar \beta, U_1^-, W)}}{W P(\bar \beta,W)^3} \Big( a e^{\int_{U_m}^U \frac{f(W)}{W^2} dW} + \delta P_2(W) \Big)^2 dW \\
& \quad + \int_{U_m}^{U^+_1} \frac{2 f(W) e^{-K(\bar \beta, U_1^-, W)}}{W P(\bar \beta,W)^3} \Big( b e^{\int_{U_m}^U \frac{f(W)}{W^2} dW} + \delta P_3(W) \Big)^2 dW.
\end{split}
\end{equation*}
Observing that from the necessary conditions for optimality \eqref{Equa:kernel_0}
\begin{equation*}
K(\bar \beta, U_1^-, U_m) = K(\bar \beta, U_m, U^+_1) = 0,
\end{equation*}
we obtain the expression
\begin{equation*}
\begin{split}
\partial^2_\beta E(\bar \beta) &= \int_{U^-_1}^{U_m} \frac{2 f(W) e^{- K(\bar \beta,U_1^-,W)}}{W P(\bar \beta,W)^3} \delta P_2(W)^2 dW + 2 a \int_{U^-_1}^{U_m} \frac{2 f(W) e^{- \int_{U_m}^{W} \frac{V f - P^2}{V P^2} dV}}{W P(\bar \beta,W)^3} e^{\int_{U_m}^W \frac{f}{V^2} dV} \delta P_2(W) dW \\
& \quad + a^2 \int_{U^-_1}^{U_m} \frac{2 f(W) e^{- \int_{U_m}^{W} \frac{V f - P^2}{V P^2} dV}}{W P(\bar \beta,W)^3} e^{2 \int_{U_m}^W \frac{f}{V^2} dV} dW + \int_{U_m}^{U^+_1} \frac{2 f(W) e^{-K(\bar \beta,U_m,W)}}{W P(\bar \beta,W)^3} \delta P_3(W)^2 dW \\
& \quad + 2 \int_{U_m}^{U^+_1} \frac{2 f(W) e^{- \int_{U_m}^{W} \frac{V f - P^2}{V P^2} dV}}{W P(\bar \beta,W)^3} b e^{\int_{U_m}^W \frac{f}{V^2} dV} \delta P_3(W) dW + \int_{U_m}^{U^+_1} \frac{2 f(W) e^{- \int_{U_m}^{W} \frac{V f - P^2}{V P^2} dV}}{W P(\bar \beta,W)^3} b^2 e^{2 \int_{U_m}^W \frac{f}{V^2} dV} dW.
\end{split}
\end{equation*}
We compute now the various terms: %using \eqref{Equa:KOIbeta}
\begin{equation*}
\begin{split}
& 2 a \int_{U^-_1}^{U_m} \frac{2 f(W) e^{- \int_{U_m}^{W} \frac{V f - P^2}{V P^2} dV}}{W P(\bar \beta,W)^3} e^{\int_{U_m}^W \frac{f}{V^2} dV} \delta P_2(W) dW \\
&= 2 a \int_{U^-_1}^{U_m} \frac{2 f(W) e^{\int_{U_m}^{W} \frac{1}{V} dV}}{W P(\bar \beta,W)^3} \delta P_2(W) dW = \frac{2 a}{U_m} \int_{U^-_1}^{U_m} \frac{2 f(W)}{W P(\bar \beta,W)^3} \delta P_2(W) dW = 0,
\end{split}
\end{equation*}
\begin{equation*}
\begin{split}
&a^2 \int_{U^-_1}^{U_m} \frac{2 f(W) e^{- \int_{U_m}^{W} \frac{V f - P^2}{V P^2} dV}}{W P(\bar \beta,W)^3} e^{2 \int_{U_m}^W \frac{f}{V^2} dV} dW = \frac{a^2}{U_m} \int_{U^-_1}^{U_m} \frac{2 f(W)}{P(\bar \beta,W)^3} e^{\int_{U_m}^W \frac{f}{V^2} dV} dW,
\end{split}
\end{equation*}
\begin{equation*}
\begin{split}
&2 \int_{U_m}^{U^+_1} \frac{2 f(W) e^{- \int_{U_m}^{W} \frac{V f - P^2}{V P^2} dV}}{W P(\bar \beta,W)^3} b e^{\int_{U_m}^W \frac{f}{V^2} dV} \delta P_3(W) dW = 2 \int_{U_m}^{U^+_1} \frac{2 f(W)}{W P(\bar \beta,W)^3} b \frac{W}{U_m} \delta P_3(W) dW = 0,
\end{split}
\end{equation*}
\begin{equation*}
\begin{split}
&\int_{U_m}^{U^+_1} \frac{2 f(W) e^{- \int_{U_m}^{W} \frac{V f - P^2}{V P^2} dV}}{W P(\bar \beta,W)^3} b^2 e^{2 \int_{U_m}^W \frac{f}{V^2} dV} dW = \frac{b^2}{U_m} \int_{U_m}^{U^+_1} \frac{2 f(W)}{P(\bar \beta,W)^3} e^{\int_{U_m}^W \frac{f(V)}{V^2} dV} dW.
\end{split}
\end{equation*}
Hence we obtain
\begin{equation*}
\begin{split}
\partial^2_\beta E(\bar \beta-) &= \int_{U^-_1}^{U_m} \frac{2 f(W) e^{-K(\bar \beta,U_1^-, W)}}{W P(\bar \beta,W)^3} \delta P_2(W)^2 dW + \int_{U_m}^{U^+_1} \frac{2 f(W) e^{-K(\bar \beta, U_m, W)}}{W P(\bar \beta,W)^3} \delta P_3(W)^2 dW \\
& \quad + \frac{a^2}{U_m} \int_{U^-_1}^{U_m} \frac{2 f(W)}{P(\bar \beta,W)^3} e^{\int_{U_m}^W \frac{f}{V^2} dV} dW + \frac{b^2}{U_m} \int_{U_m}^{U^+_1} \frac{2 f(W)}{P(\bar \beta,W)^3} e^{\int_{U_m}^W \frac{f(V)}{V^2} dV} dW.
\end{split}
\end{equation*}
We thus obtain substituting the expressions of $a,b$
\begin{equation*}
\begin{split}
\partial^2_\beta E(\bar \beta-) - \partial^2_\beta E(\bar \beta+) &= \frac{a^2}{U_m} \int_{U^-_1}^{U_m} \frac{2 f(W)}{P(\bar \beta,W)^3} e^{\int_{U^-_1}^W \frac{f}{V^2} dV} dW + \frac{b^2}{U_m} \int_{U_m}^{U^+_1} \frac{2 f(W)}{P(\bar \beta,W)^3} e^{\int_{U^+_1}^W \frac{f(V)}{V^2} dV} dW \\
&= \frac{ \big( \int_{U_m}^{U^+_1} \frac{2 f(W)}{P^3} e^{\int_{U_m}^W \frac{f}{P^2} dV} dW \big) \big( \int_{U^-_1}^{U_m} \frac{2 f(W)}{P^3} e^{\int_{U_m}^W \frac{f}{P^2} dV} dW \big)}{U_m \int_{U^-_1}^{U^+_1} \frac{2 f(W)}{P^3} e^{\int_{U_m}^W \frac{f}{P^2} dV} dW} \big( \delta P_3(U_m) - \delta P_2(U_m) \big)^2 > 0.
\end{split}
\end{equation*}
Hence there is a discontinuity in $\partial^2_\beta E(\beta)$.

\newpage

%\printglossaries

\newpage

\printbibliography

\end{document}